# 2. A comparison of two methods applied to the optimisation of fluid power circuits

A.M.Connor & D.G.Tilley

**Abstract**

*This paper describes two optimisation methods which can be applied to the parameter selection stage of Fluid Power System design. These two methods used are a Genetic Algorithm (GA) and a Tabu Search method, both of which have been claimed to be truely global methods. GAs are a method inspired by natural selection and Darwinian evolution whilst Tabu Search is an aggressive search metaheuristic which guides local search methods towards the globally optimum solution. Results are presented for two different circuit optimisation tasks. These results show that each of the two methods have both advantages and disadvantages.*

## 1    Introduction

The manual design of Fluid Power circuits is a lengthy process which often involves both a degree of trial and error, and considerable expert knowledge, in the selection and sizing of components. This often leads to a compromise in performance in order to achieve results.

The use of numerical optimisation techniques linked with dynamic simulation can lead to better solutions but often require a large amount of time to find the best solution. A significant portion of this time is due to the complex numerical integration required to carry out dynamic simulation of hydraulic circuits.

This paper describes the underlying methodology of two optimisation methods and compares the performance of each on a number of different Fluid Power circuits.

## 2 Tabu Search

Tabu Search is a metaheuristic which is used to guide local search methods towards a globally optimum solutions. The power of the Tabu Search algorithm is derived from the use of flexible memory cycles of differing time spans. These memory cycles are used to control, intensify and diversify the search in order to find a suitable solution.

### 2.1 Short Term Memory

The most simple implementation of a Tabu search is based around the use of a hill climbing algorithm. Once the method has located a locally optimal solution, the short term memory is used to force the search out from this optima in a different direction.

The short term memory constitutes a form of aggressive search that seeks to always make the best allowable move from any given position. Short term memory is implemented in the form of *tabu restrictions* which prevent the search from cycling around a given optimal or sub-optimal position. Essentially, this restriction applies to a limited set of previously visited positions which is continuously updated on a "first in, first out" basis as the search progresses.

The short term memory cycle has been implemented as a list of Tabu restrictions. The Tabu list contains the parameter values of the last 'n' accepted solutions. During the search, before the objective function of a proposed point is evaluated the parameter set is compared to those contained in the Tabu list. If the position has been recently visited then the objective function is not evaluated and the search classes the move as Tabu.

### 2.2 Search Intensification

Any intensification of the search is based around the concept of the intermediate memory cycle. In this case, the intermediate memory is very similar to the short term memory in that it is a list of previously accepted solutions. However, the solutions contained in the intermediate memory list are the previous 'm' best solutions. As a new best solution is found, this is placed into the list and the oldest solution removed. This is termed intermediate memory as the time scale of replacement is much longer than for short term memory.

## 2.3 Search Diversification

In many Tabu search applications, search diversification is achieved by the use of long term memory cycles or just by randomly refreshing the current base point when certain conditions apply.

However, in previous work [1] relating to hydraulic circuit optimisation certain parallels have been drawn to the notion of schema processing within Genetic Algorithms to propose a novel diversification strategy that is intended to force the search towards new areas of the search domain which also have good objective function values. The search is refreshed by randomly selecting values for each parameter of the search from the solutions contained in the intermediate term memory list. This approach is particularly suitable where parameters are scaled to within specified limits.

## 2.4 Hill Climbing Algorithm

The underlying hill climbing algorithm used in this work is based upon the method developed by Hooke and Jeeves [2]. This method consists of two stages. The first stage carries out an initial exploration around a given base point. Once this exploration has been carried out and a new point found, the search is extended along the same vector by a factor *k*. This is known as a pattern move. If the pattern move locates a solution with a better objective value than the exploration point, then this point is used as a new base point and the search is repeated. Otherwise, the search is repeated using the exploration point as a new base point.

In the Hooke and Jeeves search, this process is repeated until no further improvement is found. When this occurs, the step size is reduced and the search continued. This is repeated until the step size falls below a given value.

In the Tabu search, this search strategy is somewhat modified. Firstly, the selection of the new base point is adjusted so that the search always selects the best available move. When no improvement can be found, the best available move is defined as that where the increase in objective function value is the smallest. Because of this, a new control algorithm has been developed which reduces the step size as the search progresses and eventually terminates the search.

## 3   Genetic Algorithms

Genetic Algorithms [3] are a non-derivative based optimisation method based on the Darwinian theory of natural selection and survival of the fittest. The method starts with a initial population of randomly generated solutions and each successive generation is formed by mimicking the genetic operators that occur in natural systems. These operators are reproduction, crossover and mutation. As the number of generations increases the population "evolves" towards the global optimum solution.

Each individual in the population consists of a binary string which can be decoded to find the actual parameter values for the solution. This binary representation can be shown to be the most effective representation by considering how characteristics of individuals are propagated from one generation to the next. This propagation is described in the schema growth equation [4].

The main advantage of Genetic Algorithms is that they utilise a population of solutions which inhabit the solution space defined by the objective function. This inherently parallel search has the potential to explore different areas of the solution space at the same time and extract the best solution to a given problem.

The Genetic Algorithm used in this work has previously been applied to problems in the area of fluid power system design [5,6]. It is a more complex method than the traditional Genetic Algorithm as it utilises a number of parallel sub-populations which exchange individuals through a process known as migration. GAs of this type are known as Parallel Genetic Algorithms (PGA) and they have been shown to be more effective than simple GAs on complex problems.

The PGA used in this work has eight sub-populations and every three generations, four individuals migrate from each sub-population. Each sub-population also has different probabilities of crossover and mutation so as to maximise the potential for diversification with in the search.

## 4    Simple Hydrostatic Transmission

Figure 4.1 shows the circuit diagram for a simple hydrostatic transmission. The objective of the optimisation is to select capacities for the pump and motor in the hydrostatic transmission so that a given rotational speed is produced by the motor within the first four seconds.

All of the other parameters in the circuit simulation remain constant. The prime mover speed is 1500 rpm, the relief valve cracking pressure is 100 bar and the moment of inertia is 50 kgm$^2$.

Brief details of the characteristics of each of the component models used in this circuit are given in Appendix A. Hoever, further details can be found in [7].

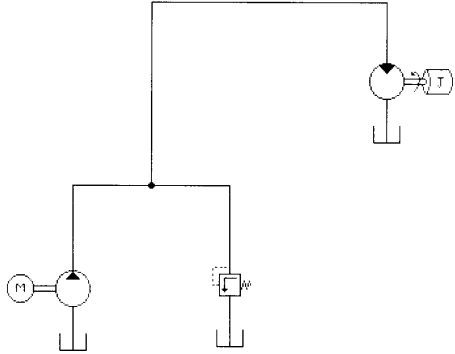

Fig 4.1 : Simple Hydrostatic Transmission

The parameters selected for inclusion in the optimisation process are the pump displacement, the motor displacement and the diameter of the connecting pipes. This differs from previous work [1] which did not include the pipe diameter as a design variable.

The objective function was calculated by using the following equation;

$$obfn = \left(\omega_{desired} - \omega_{actual}\right)^2 \times \left(1 + \frac{P_D}{P_{UB}}\right) \times \left(1 + \frac{RV_{flow}}{P_{flow}}\right)$$

where;

$\omega_{desired}$ = Desired Speed

$\omega_{actual}$ = Actual Speed

$P_D$ = Pump Displacement

$P_{UB}$ = Upper Bound on Pump Size

$P_{flow}$ = Pump Flow

$RV_{flow}$ = Relief Valve Flow

This objective function differes from that used in previous work [5,6] in which the different objective function terms were combined into a single objective function value by the use of a weighted sum approach. However, it is clear that the primary objective is the motor output speed and that the pump size and relief valve flow terms are only being used to penalise inefficient solutions. Solutions which have a small motor size are likely to have a lower power consumption for a given prime mover speed. Similarly, solution which have very small relief valve flows in comparison to the total pump flow are less likely to suffer from performance degreadation due to an increased fluis temperature.

Multiplying the prime design objective by scaled penalty functions has several advantages over a weighted sum approach. Firstly, it is uneccessary to select values for the weightings associated with the objective function components. These are normally found using a degree of trial and error whilst observing the progress of the optimisation method over a number of trial runs. The second advantage is that as the error on the prime design objective is reduced so is the significance of the penalty terms.

As both the Tabu Search and the Genetic Algorithm include a degree of probability, ten runs of each algorithm were performed on each problem with different, randomly generated starting conditions. The following solutions were found by the Tabu Search method.

| RUN | PUMP DISPLACEMENT (CC/REV) | MOTOR DISPLACEMENT (CC/REV) | PIPE DIAMETER (MM) | OBFN | EVALS |
|---|---|---|---|---|---|
| 1 | 28 | 139 | 35.5 | 0.000467 | 818 |
| 2 | 77 | 384 | 7 | 0.066865 | 709 |
| 3 | 47 | 234 | 31.5 | 0.065513 | 888 |
| 4 | 57 | 284 | 52 | 0.004882 | 921 |
| 5 | 152 | 158 | 11.5 | 0.026517 | 777 |
| 6 | 65 | 324 | 55 | 0.000082 | 701 |
| 7 | 153 | 763 | 57 | 0.023203 | 907 |
| 8 | 83 | 414 | 59 | 0.002309 | 707 |
| 9 | 185 | 923 | 46 | 0.00696 | 679 |
| 10 | 65 | 324 | 55 | 0.000082 | 877 |
| Avg | 91.2 | 394.7 | 40.95 | 0.019688 | 798.4 |
| SD | 52.78215 | 254.8895 | 19.04592 | 0.02626 | 95.43025 |

Fig 4.2 : Tabu Search Solutions

For this simple circui is is relatively simple to check he validity of each of the solutions. Assuming there are no losses in the system the relationship between motor displacement and pump displacement is given by the following equation:

$$\omega_{PM} \times P_D = \omega_P \times M_D$$

where;

$\omega_{PM}$ = Prime Mover Speed
$\omega_D$ = Motor Actual Speed
$P_D$ = Pump Displacement
$M_D$ = Motor Displacement

Given that the prime mover speed remains constant at 1500 rpm and the nominal operating speed of each of the solutions is 300 rpm it follows that the ratio of motor displacement to pump displacement should be 5:1. Examining the solutions n Figure 4.2 shows this to be the case.

The best of these solutions, which required 877 evaluations, has the speed response shown in Figure 4.3.

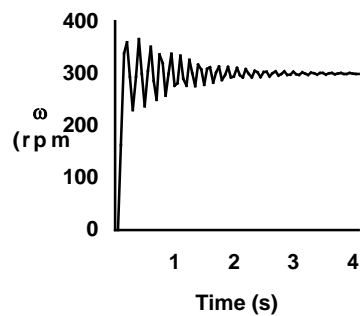

Fig 4.3 : Speed Response of Best Solution

This solution suffers from a very oscillatory transient response during start up. In this work this oscillation is not penalised as the design objective concerning the motor speed considers the value at the end of the time period and does not even require the solutions to have reached a steady state value. However, the transients die away and at steady state the motor is rotating with a nominal constant velocity of 299 rpm. This solution also suffers from another drawback. During the transient phase, the relief valve opens and shuts in a short space of time. This has not been penalised by the inclusion of the relief valve flow into the objective function as this is only considering values of state variables at the end of the simulation time. The objective function could be enhanced to consider the values of state variables throughout the entire simulation time so as to force the method towards solutions which do not exhibit this characteristic.

The following solutions were found by the Genetic Algorithm

| RUN | PUMP DISPLACEMENT (CC/REV) | MOTOR DISPLACEMENT (CC/REV) | PIPE DIAMETER (MM) | OBFN | EVALS |
|---|---|---|---|---|---|
| 1 | 193 | 963 | 46 | 0.006333 | 6882 |
| 2 | 197 | 982 | 58.5 | 0.157148 | 6943 |
| 3 | 197 | 982 | 58.5 | 0.157148 | 7018 |
| 4 | 193 | 963 | 37 | 0.006018 | 6916 |
| 5 | 197 | 982 | 58.5 | 0.157148 | 7028 |
| 6 | 193 | 963 | 40 | 0.005727 | 6923 |
| 7 | 189 | 943 | 60 | 0.007789 | 6948 |
| 8 | 197 | 982 | 58.5 | 0.157148 | 6917 |
| 9 | 197 | 982 | 58.5 | 0.157148 | 6936 |
| 10 | 193 | 963 | 9 | 0.005657 | 6985 |
| Avg | 194.6 | 970.5 | 48.45 | 0.081726 | 6949.6 |
| SD | 2.796824 | 13.50926 | 16.3409 | 0.079504 | 46.87857 |

Fig 4.4 : Genetic Algorithm Solutions

Again, each of these solutions exhibits a motor displacement to pump displacement ration of approximately 5:1. However, the GA appears to be locating solutions which are clustered within a much smaller region of the solution space.

The best solution, which required 6985 evaluations, has the speed response shown in Figure 4.5.

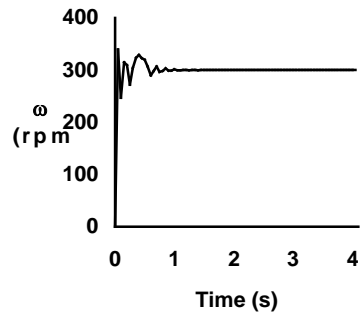

Fig 4.5 : Speed Response of Best Solution

This solution reaches a steady state operating speed of 300rpm in a much shorter space of time than the solution found by the Tabu Search. Also, in this solution the relief valve does not open at all during the simulation. Comparing the performance of the methods, it can be seen that the Tabu Search method is finding solutions, on average, with lower objective function values than the GA. Yet examining the performance of the best solutions found by each method it can be seen that the GA has located a solution that offers the best practical characteristics. This is caused by a lack of definition in the objective function.

By considering the standard deviation of the parameter values, it can also be seen that the GA is more consistent than the Tabu Search method. However, by considering the actual values of the parameters it can be seen that the problem representation could be improved. In general both methods are finding solutions where the pipe diameters are quite large. Reducing the allowed boundaries on certain variables should force the methods towards more appropriate solutions.

## 5     Regenerative Hydrostatic Transmission

Figure 5.1 shows an alternative closed loop hydrostatic transmission which includes a regenerative feed circuit. In this optimisation, the circuit has been modelled so that the motor is suffering from high slip loss which represents a faulty condition. The motor displacement is not selected as a design variable so that the volumetric and mechanical efficiencies are not changed from the values defined by specifying the slip loss parameters. Again. The prime mover speeds are 1500 rpm, the relief valvle cracking pressures are 100 bar and the load inertia is 50 kgm$^2$. Details of the components models are in Appendix A.

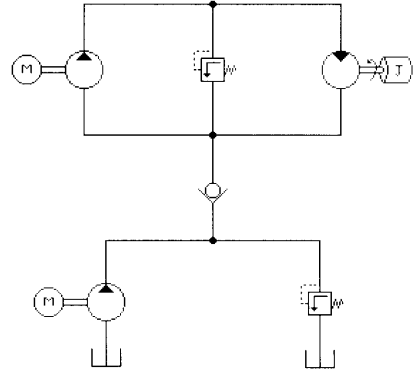

Fig 5.1 : Closed Loop Transmission with Boost Pump

The aim of this optimisation is to select displacements for both of the pumps in the circuit as well as the operating speed of the associated prime movers. The design variables are being found when a fault condition exists in the motor. The optimisation is carried out so that even when a fault occurs the operating speed of the motor does not fall below a given threshold. This threshold has been modelled as a volumetric efficiency of approximately 75%. Optimised solutions can then be compared for when the fault does and does not exist. The objective function was calculated by using the following equation;

$$obfn = (\omega_{desired} - \omega_{actual})^2 \times \left(1 + \frac{RV_{flow}}{Pf_{low}}\right)_{main} \times \left(1 + \frac{RV_{flow}}{Pf_{low}}\right)_{feeder}$$

where;

$\omega_{desired}$ = Desired Speed

$\omega_{actual}$ = Actual Speed

$P_{flow}$ = Pump Flow

$RV_{flow}$ = Relief Valve Flow

The nature of this objective function is similar to that used in the previous example. The main different is that the two penalty function terms are based around the flow through the relief value in the circuit. In this example, the size of the two pumps are not included although future work could consider the addition of two more penalty terms based on the sizes of the two motore.

21

Solutions found by the Tabu Search method are shown in Figure 5.2.

| RUN | PUMP 1 DISPLACEMENT (CC/REV) | $W_1$ (RPM) | PUMP 2 DISPLACEMENT (CC/REV) | $W_2$ (RPM) | OBFN | EVALS |
|---|---|---|---|---|---|---|
| 1 | 10 | 1221 | 143 | 1327 | 2.17E-05 | 1197 |
| 2 | 136 | 164 | 394 | 484 | 2.90E-06 | 1079 |
| 3 | 39 | 730 | 398 | 480 | 0.000244 | 1114 |
| 4 | 26 | 1029 | 326 | 585 | 6.88E-07 | 1195 |
| 5 | 43 | 678 | 696 | 276 | 1.47E-07 | 1570 |
| 6 | 95 | 265 | 138 | 1377 | 2.54E-07 | 1204 |
| 7 | 157 | 105 | 163 | 1165 | 4.06E-06 | 1097 |
| 8 | 128 | 228 | 475 | 403 | 0.043869 | 1332 |
| 9 | 249 | 100 | 429 | 445 | 2.05E-06 | 1282 |
| 10 | 157 | 105 | 163 | 1165 | 4.06E-06 | 1097 |
| Avg | 104 | 462.5 | 332.5 | 770.7 | 0.004415 | 1216.7 |
| SD | 75.41736 | 419.571 | 182.8984 | 431.4335 | 0.013863 | 149.653 |

Fig 5.2 : Tabu Search Solutions

The best solution, which required 1570 evaluations, has the speed response shown in Figure 5.3 for the faulty and in Figure 5.4 for the non-faulty condition.

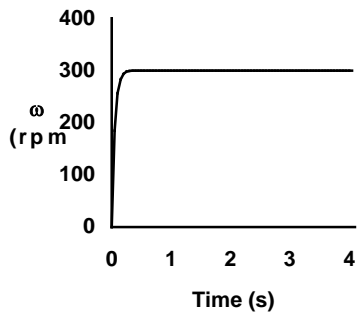 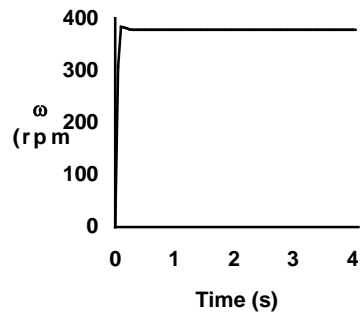

Fig 5.3 : Speed Response with Fault        Fig 5.4 : Speed Response without Fault

When the fault is being modelled, the optimal solution has a steady state operating speed of 300rpm. The flow from the boost pump is 28.9 L/min and main pump flow is 189.4 L/min. Neither of the relief valves open.

For when the motor is being modelled so as to not include the fault, the steady state operating speed is 378rpm. The flow from the boost pump is 28.9 L/min, of which 24.7 L/min returns directly to the tank through the relief valve. The main pump flow is 189.3 L/min and the corresponding relief valve opens briefly during the transient start up phase.

The solutions found by the Genetic Algorithm are shown in Figure 5.5.

| RUN | PUMP 1 DISPLACEMENT (CC/REV) | $W_1$ (RPM) | PUMP 2 DISPLACEMENT (CC/REV) | $W_2$ (RPM) | OBFN | EVALS |
|---|---|---|---|---|---|---|
| 1 | 275 | 104 | 227 | 837 | 0.530079 | 7404 |
| 2 | 103 | 115 | 352 | 541 | 0.329347 | 7470 |
| 3 | 10 | 139 | 282 | 673 | 0.014378 | 7267 |
| 4 | 270 | 100 | 196 | 969 | 0.346373 | 7353 |
| 5 | 46 | 152 | 200 | 947 | 0.343346 | 7141 |
| 6 | 10 | 152 | 282 | 673 | 0.014356 | 7185 |
| 7 | 213 | 100 | 103 | 1839 | 0.599795 | 7348 |
| 8 | 206 | 125 | 363 | 525 | 0.105261 | 7419 |
| 9 | 46 | 152 | 200 | 947 | 0.343346 | 7330 |
| 10 | 184 | 100 | 328 | 580 | 0.069611 | 7148 |
| Avg | 136.3 | 123.9 | 253.3 | 853.1 | 0.269589 | 7306.5 |
| SD | 105.119 | 23.06248 | 82.49586 | 386.9879 | 0.209367 | 116.5077 |

Fig 5.5 : Genetic Algorithm Solutions

In this example, the GA is not as consistent as in the first example. This is probably due to a poor objective function definition.

The best solution, which required 7185 evaluations, has the speed shown in Figure 5.6 for the faulty and in Figure 5.7 for the non-faulty condition.

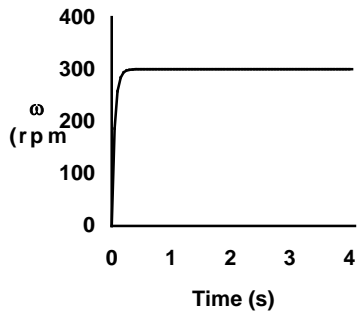 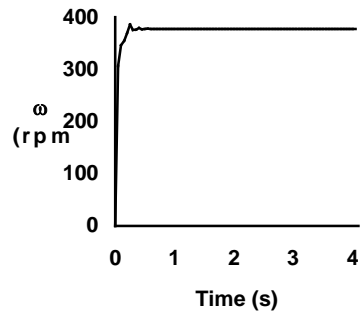

Fig 5.6 : Speed Response with Fault       Fig 5.7 : Speed Response without Fault

When the fault is being modelled, the optimal solution has a steady state operating speed of 300.1rpm. The flow from the boost pump is 1.4 L/min and main pump flow is 189.5 L/min. The required flow from the boost pump is very low. However, neither of the relief valves open.

For when the motor is being modelled so as to not include the fault, the steady state operating speed is 378rpm.The flow from the boost pump is 1.4 L/min, which again is very low. The main pump flow is 189.0 L/min and the corresponding relief valve opens briefly during the transient start up phase.

The best solution found by the Tabu Search has a lower objective function value than that of the best solution found by the GA. However, by considering the performance of each circuit under the two extreme modelled conditions it can be seen that the solution found by the GA is likely to be a more fault tolerant system as the flow through the relief valves are much lower when a fault is not being modelled. However, the very low required flow from the boost pump is unlikely to be obtainable and so this solution may not be practically implementable. The solution found by the Tabu Search has quite high relief valve flows which will lead to an increased fluid temperature and a corresponding degradation in performance. Neither solution is then ideal and more work is required in defining design objective into a form which can be utilised by a numerical optimisation technique in an effective manner..

This type of optimisation task may be carried out more effectively by a better definition of the design objectives. The use of fuzzy logic may lead to the possibility of defining objectives for multiple operating conditions in such a way as to make the solutions discovered by either method more robust to changes in operating condition.

## 6   Discussion

The results presented in the previous two sections show that both of the methods described in this paper are suitable for application to the optimisation of Fluid Power systems. Both methods are suffering slightly from poor objective function definition, though in all cases both methods are locating solutions which perform satisfactorily. Each method has both advantages and disadvantages. By comparing these, it can be seen that the methods could be used to complement each other in a automatic Fluid Power system design methodology.

The great disadvantage of the PGA is the large number of evaluations required to find a solution. This is due to the lack of termination criteria built into the method. The PGA runs for a fixed number of generations and then performs a local search in the region of the best found solution. The inclusion of termination criteria in a GA has been shown [7] to bias the action of the method and often leads to the location of sub-optimal solutions.

The Tabu Search, however, is an aggressive search which requires less objective function evaluations to find a solution. It is also more suitable for discrete parameter optimisation than the GA and it is straight forward to build in termination criteria into the search with out biasing the final solution.

## 7  Conclusions

The results presented in this paper indicate that both of methods considered have relative merits. Tabu Search offers an aggressive search that locates acceptable solutions in a very small number of objective function evaluations. However, the GA tends to locate solutions in a more consistent manner.

Further work is being carried out in the following areas. The first is to develop tools which allow the design objectives to be formulated into well balanced objective functions. This is likely to involve a degree of intelligence embedded into a pre-processor for optimisation of Fluid Power systems.

The Tabu Search method appears to offer the most efficient optimisation of Fluid Power systems when coupled with a well defined objective function. Work will continue to enhance the method and make it more robust.

## 8  Acknowledgments

The research in this paper is funded by the Engineering and Physical Sciences Research Council under grant no. 86796. This support is gratefully acknowledged.

**Appendix A**

This appendix provides details of the component models used in the simulation and optimization of the circuits described. The same models are used for components common to both circuits. Each component is given a model number taken from [7] for future reference.

The pump model (PU01) is an instantaneous model of a fixed displacement hydraulic pump which may be gear, vane or piston variety. The model considers the steady state behavior of the pump and assumes that the flow is compressible. Losses are lumped into three terms, a slip loss (leakage to drain and cross port follow), a pressure dependent friction loss and a speed dependent friction loss.

The motor model (MO01) is a companion model to the PU01 model. Again the flow is assumed to be incompressible and the model uses the same losses terminology.

The relief valve model (RV00) can be used to model either single stage or a two stage relief valves. The flow is assumed to be incompressible and the model does not take into account the dynamics of the valve. In addition to this is is assumed that pressure rise varies linearly to the flow rate once the cracking pressure has been exceeded.

The prime mover model (PM01) can be used to model either an electrical or diesel prime mover connected to a pump. The model does not include any transient calculation and is assumed to apply a constant speed.

The pipe model (HP00) is a dynamic, frictionless, constrant volume, compressible hydraulic pipe model which includes models of both cavitation and air release. Pipe friction and fluid intertia effects are ignored. The fluid compressibility is modelled as a first order system.

The load model (RL01) is a model of a rotary load which incorporates stiction, Coulomb friction, speed dependent friction and windage. There is also provision for a constant applied load torque. The effects of temperature and ageing are not taken into account.

The check valve model (CV00) is a model of a hydraulic free or spring operated check valve. The valve is assumed to operate instantaneously and has a liner flow rate/pressure relationship with the opening point modelled as a discontinuity.